\documentclass[reqno,11pt]{amsart} 
\usepackage{amsmath}
\usepackage{amssymb}
\usepackage{amsthm}
\usepackage{verbatim}
\usepackage{xcolor}
\usepackage{graphics}

\newcommand\NoBlackBoxes{\global\overfullrule0pt}

\NoBlackBoxes
\theoremstyle{plain}

\begin{document}

\begin{abstract}
Berry-Esseen-type bounds are developed in the multidimensional
local limit theorem in terms of the Lyapunov coefficients and
maxima of involved densities.
\end{abstract}

\title[Local Limit Theorems]{Berry-Esseen bounds \\
in local limit theorems}

\author{S. G. Bobkov$^{1,4}$}
\thanks{1) 
School of Mathematics, University of Minnesota, Minneapolis, MN, USA,
bobkov@math.umn.edu. 
}

\author{F. G\"otze$^{2,4}$}
\thanks{2) Faculty of Mathematics,
Bielefeld University, Germany,
goetze@math-uni.bielefeld.de.
}

\thanks{4) Research supported by the NSF grant DMS-2154001 and 
the GRF – SFB 1283/2 2021 – 317210226 }


\maketitle
\markboth{S. G. Bobkov and F. G\"otze}{Berry-Esseen bounds}

\def\theequation{\thesection.\arabic{equation}}
\def\E{{\mathbb E}}
\def\R{{\mathbb R}}
\def\C{{\mathbb C}}
\def\P{{\mathbb P}}
\def\Z{{\mathbb Z}}
\def\S{{\mathbb S}}
\def\I{{\mathbb I}}
\def\T{{\mathbb T}}

\def\s{{\mathbb s}}

\def\G{\Gamma}

\def\Ent{{\rm Ent}}
\def\var{{\rm Var}}
\def\Var{{\rm Var}}
\def\V{{\rm V}}

\def\H{{\rm H}}
\def\Im{{\rm Im}}
\def\Tr{{\rm Tr}}
\def\s{{\mathfrak s}}

\def\k{{\kappa}}
\def\M{{\cal M}}
\def\Var{{\rm Var}}
\def\Ent{{\rm Ent}}
\def\O{{\rm Osc}_\mu}

\def\ep{\varepsilon}
\def\phi{\varphi}
\def\vp{\varphi}
\def\F{{\cal F}}

\def\be{\begin{equation}}
\def\en{\end{equation}}
\def\bee{\begin{eqnarray*}}
\def\ene{\end{eqnarray*}}

\thispagestyle{empty}

\section{{\bf Introduction}}
\setcounter{equation}{0}

\vskip2mm
\noindent
Consider the normalized sum
$$
Z_n = \frac{X_1 + \dots + X_n}{\sqrt{n}}
$$
of $n$ independent random vectors $X_k$ in $\R^d$ with mean zero and
covariance matrices $\sigma_k^2 I_d$ ($\sigma_k > 0$) such that 
$\sigma_1^2 + \dots + \sigma_n^2 = n$. The latter ensures that 
$Z_n$ has mean zero and unit covariance matrix $I_d$.

Introduce the Lyapunov ratios of the third and fourth orders
$$
\beta_3 = \sup_{|\theta| = 1}\,
\bigg[\,\frac{1}{n} \sum_{k=1}^{n} \E\,|\left<\theta,X_k\right>|^3\bigg], 
\quad
\beta_4 = \sup_{|\theta| = 1}\,
\bigg[\,\frac{1}{n} \sum_{k=1}^{n} \E \left<\theta,X_k\right>^4\bigg],
$$
assuming respectively that $\E\,|X_k|^3 < \infty$ and $\E\,|X_k|^4 < \infty$
for all $k \leq n$. In dimension one, the quantity $\beta_3$ is commonly 
used to quantify the normal approximation for the distribution of $Z_n$ 
in a weak sense via the Berry-Esseen bound
$$
\sup_{x \in \R}\, \big|\P\{Z_n \leq x\} - \P\{Z \leq x\}\big| \leq 
\frac{C}{\sqrt{n}}\,\beta_3,
$$
where $Z$ is a standard normal random variable and $C$ is an absolute
constant (cf. e.g. \cite{P}). The Lyapunov ratio $\beta_4$ also appears 
in a natural way, for example, for the approximation of the characteristic function 
of $Z_n$ by a corrected normal characteristic function. Both quantities may also 
be used to control the distance between the distribution of $Z_n$ and $Z$ on 
the real line in total variation and in relative entropy, cf. \cite{B-C-G2}.

The aim of this note is to quantify the normal approximation in a stronger
sense towards a uniform 
convergence of densities $p_n$ of $Z_n$ to the standard normal density
$$
\varphi(x) = \frac{1}{(2\pi)^{d/2}}\,e^{-|x|^2/2}, \quad x \in \R^d.
$$
In the i.i.d. situation (when all $X_k$ are identically distributed with 
$\sigma_k=1$), the necessary and sufficient condition for the
convergence of the uniform distance
$$
\Delta_n = \sup_x |p_n(x) - \varphi(x)|
$$
to zero as $n \rightarrow \infty$ is that $p_n$ is bounded for some 
$n = n_0$ (cf. \cite{P,B-RR}). As typically $n_0=1$ 
in applications, it is natural to consider the case where all $X_k$ have 
bounded densities. 

So, introduce the maximum-of-density functional
$M(X) = {\rm ess\,sup}_x \, p(x)$,
where $p$ denotes a density of a random vector $X$, and define
$$
M = \max_k M(X_k), \quad \sigma^2 = \max_k \sigma_k^2 \ \ (\sigma>0).
$$
The functionals $\beta_3$, $\beta_4$, $\sigma$ and $M$ can be used 
to derive the following upper bounds (which seem to be new already in the
one dimensional situation).

\vskip5mm
{\bf Theorem 1.1.} {\sl With some positive absolute constant $C$, the density
$p_n$ of $Z_n$ satisfies
\be
\Delta_n \leq \frac{(C \sigma)^d M^2}{\sqrt{n}}\,\beta_3.
\en
Moreover, if $\E \left<\theta,X_k\right>^3 = 0$ for all $\theta \in \R^d$ and 
$k \leq n$ (in particular, if the distributions of $X_k$ are symmetric), then
\be
\Delta_n \leq \frac{(C\sigma)^{2d} M^3}{n}\, \beta_4.
\en
}

\vskip2mm
Hence, when $M$ and $\sigma$ are bounded, it is possible to strengthen the Berry-Esseen
theorem with an extension to higher dimensions.

In order to reflect the influence of $M(X_k)$ on average (similarly to $\beta_3$), rather 
than via the maximal value $M$, some refinements of the bounds (1.1)-(1.2) are given 
in Section 6.

For several classes of probability distributions, the functional $M(X)$ is of the order 
$\sigma^{-d}$ when a random vector $X$ has a covariance matrix $\sigma^d I_d$
(modulo $d$-dependent constants). Here is an example involving convexity properties
of distributions.

\vskip5mm
{\bf Corollary 1.2.} {\sl  Suppose that the random vectors $X_k$ have log-concave
densities, with mean zero and unit covariance matrix. Then
\be
\Delta_n \leq \frac{C_d}{\sqrt{n}}
\en
with some constant $C_d$ depending on $d$ only. If additionally the distributions 
of $X_k$'s are symmetric, then
\be
\Delta_n \leq \frac{C_d}{n}.
\en
}

Based on the application of the Fourier transform (which is typical in the study
of various variants of the central limit theorems), the main arguments used in the proof 
of Theorem 1.1 employ the subbaditivity property of the maximum-of-density
functional $M(X)$ with respect to convolutions. This tool has been introduced in
the field of limit theorems only recently and is used in \cite{B-C-G3}, \cite{B-G}.
Another ingredient in the proof involves the extension of the Statuljavičus separation
theorem for characteristic functions to higher dimensions.

The paper is organized as follows. We start with remarks on isotropic constants
and general bounds for $M(X)$ in terms of the covariance matrix of a random vector 
$X$ (Section~2) and then discuss the question of how one can separate
the absolute value of a given characteristic function $f(t)$ from 1 in terms of $M(X)$ 
outside a neighborhood of the origin (Sections 3-4). In Sections 5 we recall the
subadditivity property of $M(X)$ and develop its applications to the integrability
properties of powers of characteristic functions. Basic results on the normal 
approximation of products of characteristic functions are recalled
in Section~6. In Section 7, we derive refined bounds on $\Delta_n$, which are
used in Section 8 to finish the proof of Theorem~1.1 and Corollary 1.2.

\vskip7mm
\section{{\bf Lower Bounds on Maximum of Density via Covariance Matrix}}
\setcounter{equation}{0}

\vskip2mm
\noindent
To start with, first let us recall the general relation
\be
M^2 \sigma^2 \geq \frac{1}{12},
\en
which holds for any random variable $X$ with standard deviation $\sigma$
and maximum of density $M = M(X)$. As an early reference one can mention 
the paper by Statuljavičus \cite{St}, p.\,651, where (2.1) is stated without derivation
as an obvious fact.
In the muldimensional situation, (2.1) is extended in the form
\be
M^{2/d} \sigma^2 \geq \frac{1}{d+2}\,\omega_d^{-2/d},
\en
assuming that the random vector $X$ has a covariance matrix $\sigma^2 I_d$
(such distributions are called isotropic).
Here an equality is attained for the uniform distribution on every Euclidean ball
in $\R^d$ and, in particular, for the unit ball $B_2 = \{x \in \R^d: |x| \leq 1\}$ 
with volume
\be
\omega_d = {\rm vol}_d(B_2) = \frac{\pi^{d/2}}{\Gamma(\frac{d}{2}+1)}.
\en

This extremal property of balls has been investigated in Convex Geometry 
in the context of bounding volume of slices of convex bodies, see 
Hensley \cite{H} and Ball \cite{Ball2}. More precisely, it was shown in
\cite{Ball2}, Lemma 6, that, if the density of $X$ satisfies $p(x) \leq p(0)$
for all $x \in \R^d$, then
$$
p(0)^{2/d} \int_{\R^d} |x|^2\,dx \geq \frac{d}{d+2}\,\omega_d^{-2/d}.
$$
This amounts to (2.2) in the case where $X$ has mean zero, covariance matrix
$\sigma^2 I_d$, and assuming that $p(x)$ is maximized at the origin.
In \cite{B-M}, Proposition III.1, this inequality was generalized and 
strengthened as the property that, for any non-decreasing function 
$H = H(t)$ in $t \geq 0$, the moment-type functional
$$
\int_{\R^d} H(M^{1/d}\, |x|)\,p(x)\,dx
$$
is minimized for the uniform distribution on every Euclidean ball with center
at the origin. For $H(t) = t^2$, this property implies (2.2), assuming that 
$X$ has mean zero and covariance matrix $\sigma^2 I_d$.
But the inequality (2.2) is translation invariant, so it continues
to hold without the mean zero assumption.

The quantity $M^{1/d} \sigma$ is often 
called an isotropic constant of the distribution of $X$. It is also well-known
that the relation (2.2) is in essense dimension-free, since it implies
\be
M^{2/d} \sigma^2 \geq \frac{1}{2\pi e}
\en
with an optimal constant on the right-hand side 
(attainable asymptotically for growing dimension $d$).
More generally, if a random vector $X$ has covariance matrix $R$, then 
$$
\big(M^2\,{\rm det}(R)\big)^{1/d} \geq \frac{1}{2\pi e}.
$$
This follows from the previous lower bound (2.4) by applying it to the
random vector $Y = R^{-1/2} X$. It is isotropic with $\sigma(Y) = 1$
and $M(Y) = \sqrt{{\rm det}(R)}\, M(X)$.

In this connection, let us mention that in the class of isotropic distributions
on $\R^d$ having log-concave densities, all these lower bounds can be 
reversed as
\be
M^{2/d} \sigma^2 \leq K_d^2
\en
in terms of the maximal isotropic constant over this class for a fixed dimension $d$.
In the equivalent form, the hyperplane conjecture raised in 1980's by Bourgain, 
which is still open, asserts that $K_d$ is bounded by an absolute constant.
This is true, for example, when additionally the density of $X$ is symmetric
about all coordinate axes (cf. e.g. \cite{B-N}). As for the general log-concave case,
at this moment the best result in this direction belongs to Klartag \cite{Kl}
with his bound $K_d^2 \leq C\log(d+1)$ for some absolute constant $C$.
Let us refer an interested reader to \cite{B-G-V-V} and \cite{A-G-M} for the
history of the problem and related results.

\vskip7mm
\section{{\bf  Separation of Characteristic Functions (Statuljavičus Theorem)}}
\setcounter{equation}{0}

\vskip2mm
\noindent
In order to apply the Fourier methods for the derivation of density bounds 
in a quantitative way, one has to realize how to separate the characteristic function
$$
f(t) = \E\,e^{itX}, \quad t \in \R,
$$ 
of a random variable $X$
from 1 outside a neighborhood of the origin (which is potentially possible due
to the Riemeann-Lebesgue lemma). That is, the task is to derive estimates like
$\sup_{|t| \geq t_0} |f(t)| < 1 - \delta$ with an arbitrary $t_0 > 0$ and some positive 
$\delta = \delta(t_0)$. An important step in this direction was made by
Statuljavičus \cite{St} who derived an upper bound which we prefer to state
in an equivalent form.

\vskip5mm
{\bf Proposition 3.1.} {\sl Given a random variable $X$ with standard
deviation $\sigma$ and finite maximum of density $M = M(X)$, its characteristic 
function satisfies
\be
|f(t)| \leq 1 - \frac{c}{M^2\sigma^2}\,\min\{\sigma^2 |t|^2,1\}, \quad
t \in \R,
\en
with some absolute constant $c>0$.
}

\vskip4mm
Some generalizations of this result to the general case of unbounded densities and 
without moment assumptions are discussed in \cite{B-C-G1}. In fact,
Statulevičius considered more complicated quantities for upper bounds
reflecting the behavior of the density $p$ of $X$ on non-overlapping intervals 
of the real line and obtained the bound
$$
|f(t)| \leq 
\exp\Big\{-\frac{t^2}{96\,M^2\,(2\sigma |t| + \pi)^2} \Big\}
$$
as a partial case of a more general relation.
Therefore, let us include a shorter simplified argument aimed 
at the bound (3.1) only (without polishing the constant $c$).

\vskip5mm
{\bf Proof.}  In a more flexible form, the inequality (2.1)
can be equivalently written as
\be
\big({\rm ess\,sup}_x \, q(x)\big)^2 \int_{-\infty}^\infty x^2 q(x)\,dx
 \geq \frac{1}{12}\,\Big(\int_{-\infty}^\infty q(x)\,dx\Big)^3,
\en
holding true for any non-negative Borel measurable function $q(x)$
on the real line.

Consider the symmetrized random variable $\widetilde X = X - X'$, where 
$X'$ is an independent copy of $X$. It has a positive characteristic function 
$|f(t)|^2$ and density
$$
w(x) = \int_{-\infty}^\infty p(x+y) p(y)\,dy \leq M.
$$
Write
$$
1 - |f(2\pi t)|^2 = \int_{-\infty}^\infty (1 - \cos(2\pi tx)\, w(x)\,dx = 
2 \int_{-\infty}^\infty \sin^2(\pi tx)\, w(x)\,dx.
$$
In order to bound the last integral from below, one may use the elementary inequality
$$
|\sin(\pi x)| \geq \|x\| = \min\{|x-k|: k \in \Z\},
$$
where both sides represent 1-periodic functions. Assuming without loss
of generality that $t>0$ and using $1 - |f(s)|^2 \leq 2\,(1 - |f(s))$, this gives
\be
1 - |f(2\pi t)| \geq 4 \int_W \|tx\|^2\, w(x)\,dx,
\en
where we restrict the integration to the set
$W = \{x \in \R: t |x| < N + \frac{1}{2}\}$ for a suitable integer $N \geq 0$. 
Let us split the integral in (3.3) into the sets
$$
W_k = \Big\{x \in \R: k - \frac{1}{2} < t |x| < k + \frac{1}{2}\Big\}
$$
and rewrite (3.3) as
\bee
1 - |f(2\pi t)| 
 & \geq &
4 \sum_{k=-N}^N \int_{W_k} |tx - k|^2\, w(x)\,dx \\
 & = &
4t^2 \sum_{k=-N}^N \int_{-\frac{1}{2t}}^{\frac{1}{2t}} 
y^2\, w\Big(y + \frac{k}{t}\Big)\,dy.
\ene
Applying (3.2) to the functions
$q_k(y) = w(y + \frac{k}{t})\,1_{[-\frac{1}{2t},\frac{1}{2t}]}(y)$
and using $w \leq M$, we have
$$
\int_{-\frac{1}{2t}}^{\frac{1}{2t}} y^2 w\Big(y + \frac{k}{t}\Big)\,dy \geq
\frac{1}{12 M^2}\,
\bigg[\int_{-\frac{1}{2t}}^{\frac{1}{2t}} w\Big(y + \frac{k}{t}\Big)\,dy\bigg]^3 =
\frac{1}{12\,M^2}\, \bigg(\int_{W_k} w(x)\,dx\bigg)^3.
$$
Hence
$$
1 - |f(2\pi t)| \geq \frac{t^2}{3M^2}\, \sum_{k=-N}^N Q_k^3, \quad {\rm where} \ \ 
Q_k = \int_{W_k} w(x)\,dx.
$$
Subject to $\sum_{k=-N}^N Q_k = Q$, the sum
$\sum_{k=-N}^N Q_k^3$ is minimized when $Q_k = Q/(2N+1)$. This leads to
\be
1 - |f(2\pi t)| \geq \frac{1}{3M^2}\, \frac{t^2}{(2N+1)^2}\,Q^3.
\en

One should now maximize the right-hand side or choose a suitable $N$.
Since $\E \widetilde X = 0$, $\E \widetilde X^2 = 2\sigma^2$, we get, 
by Chebyshev's inequality,
\be
1- Q = \P\Big\{|\widetilde X| \geq \frac{N + \frac{1}{2}}{t}\Big\} \leq 2\,
\bigg(\frac{\sigma t}{N + \frac{1}{2}}\bigg)^2.
\en
If $\sigma t > \frac{1}{4}$, we choose $N = [2\sigma t + \frac{1}{2}]$, in which case
$Q \geq \frac{1}{2}$ and $2N+1 \leq 12\,\sigma t$.
Hence
$$
\frac{t^2}{(2N+1)^2}\,Q^3 \geq \frac{t^2}{(12\,\sigma t)^2 \cdot 8} =
\frac{1}{1152\,\sigma^2}
$$
and thus
$$
1 - |f(2\pi t)| \geq \frac{1}{3456\, M^2 \sigma^2}, \quad \sigma t \geq 1/4.
$$
If $\sigma t \leq \frac{1}{4}$, the choice $N = 0$ in (3.5) yields
$Q \geq \frac{1}{2}$. By (3.4),
$$
1 - |f(2\pi t)| \geq \frac{t^2}{24\,M^2}, 
\quad \sigma t \leq \frac{1}{4}.
$$
\qed

\vskip7mm
\section{{\bf Separation of Characteristic Functions (Reduction to Dimension One)}}
\setcounter{equation}{0}

\vskip2mm
\noindent
In the multidimensional case, the characteristic function
$$
f(t) = \E\,e^{i\left<t,X\right>}, \quad t \in \R^d,
$$ 
admits a similar bound.
The next statement is a preliminary step in the proof of Theorem~1.1.

\vskip5mm
{\bf Proposition 4.1.} {\sl Given a random vector $X$ in $\R^d$ with
covariance matrix $\sigma^2 I_d$, $\sigma > 0$, and a finite maximum
of density $M = M(X)$, its characteristic function satisfies
\be
|f(t)| \leq 1 - \frac{c^d}{M^2\sigma^{2d}}\,\min\{\sigma^2 |t|^2,1\}, \quad
t \in \R^d,
\en
with some absolute constant $c>0$.
}

\vskip5mm
Note that the isotropic constant $M^{1/d} \sigma$ is present on the right-hand
side, which is bounded away from zero according to (2.4).
For isotropic log-concave distributions, it is also bounded from above by a $d$-dependent
constant according to (2.5). But in this case one can certainly obtain better bounds
with decay of $f(t)$ at infinity.

For the proof, it looks natural to apply the one dimensional result to random variables 
$X_\theta = \left<\theta,X\right>$ with unit vectors $\theta$. However, it may happen 
that $X_\theta$ will have an unbounded density, which means that Proposition 3.1 
is not applicable. 

For example, in dimension $d=2$, suppose that $X = (X_1,X_2)$ has a uniform
distribution on the unbounded region $R = \{(x_1,x_2): |x_1| \leq \exp(-c|x_2|)\}$, 
$c>0$, that is, with density
$$
p(x_1,x_2) = 
\frac{c}{4}\,1_R(x_1,x_2), \quad x_1,x_2 \in \R.
$$
Then $X_1$ takes values in $(-1,1)$ and has density
$$
p_1(x_1) = \int_{-\infty}^\infty p(x_1,x_2)\,dx_2 = 
\frac{1}{2}\,\log \frac{1}{|x_1|}, \quad
-1 < x_1 < 1,
$$
which is unbounded near zero. It is easy to check that the random vector 
$X$ is isotropic, that is, $\E X_1 X_2 = 0$ and $\E X_1^2 = \E  X_2^2$ 
for a suitable constant $c$. 

\vskip5mm
{\bf Proof of Proposition 4.1.}
The above example shows that a preliminary density truncation is desirable. 
Introduce a random vector $X_r$ in $\R^d$ with parameter $r>0$ with density
$$
p_r(x) = \frac{1}{b_r}\,p_r(x)\,1_{\{|x| < r\}}, \quad x \in \R^d,
$$
where $b_r = \P\{|X| < r\}$ is a normalizing constant. Assuming that
$d \geq 2$, we choose $r = \sigma \sqrt{2d}$, which guarantees, by Markov's 
inequality, that
$$
b_r = 1 - \P\{|X| \geq r\} \, \geq \, 1 - \frac{\E\,|X|^2}{r^2} \, = \,
1 - \frac{d \sigma^2}{r^2} \, = \, \frac{1}{2}.
$$

Put $t = s\theta$, $s \in \R$, $|\theta| = 1$, and consider the characteristic functions 
$$
g_r(s) = \E\,e^{i\left<t,X_r\right>} = \E\,e^{is\left<\theta,X_r\right>}, \quad
g(s) = \E\,e^{i\left<t,X\right>} = \E\,e^{is\left<\theta,X\right>}.
$$
By construction, 
\bee
1 - |g(s)|^2 
 & = &
2 \int_{\R^d} \sin^2\Big(\frac{\left<t,x\right>}{2}\Big)\,p(x)\,dx \\
 & \geq &
2 \int_{|x| < r} \sin^2\Big(\frac{\left<t,x\right>}{2}\Big)\,p(x)\,dx \\
 & = &
2b_r \int_{\R^d} \sin^2\Big(\frac{\left<t,x\right>}{2}\Big)\,p_r(x)\,dx \\
 & \geq &
\int_{\R^d} \sin^2\Big(\frac{\left<t,x\right>}{2}\Big)\,p_r(x)\,dx \, = \,
\frac{1}{2}\,(1 - |g_r(s)|^2).
\ene
Thus,
\be
1 - |g(s)|^2 \geq \frac{1}{2}\,(1 - |g_r(s)|^2).
\en

In order to bound from below the right-hand side, first note that
\bee
\Var(\left<\theta,X_r\right>)
 & = &
\frac{1}{2b_r^2} \int_{|x|<r} \int_{|y|<r} \left<\theta,x-y\right>^2\,
p(x) p(y)\,dx dy \\
 & \leq & 
\frac{1}{2b_r^2} \int_{\R^d} \int_{\R^d} \left<\theta,x-y\right>^2\,
p(x) p(y)\,dx dy \, = \, 
\frac{\sigma^2}{b_r^2} \, \leq \, 4\sigma^2.
\ene
Thus,
\be
\Var(\left<\theta,X_r\right>) \leq 4\sigma^2.
\en

The maximum of density of the random variable $\left<\theta,X_r\right>$ 
can also be related to $M=M(X)$. For simplicity, let
$\theta = e_1 = (1,0,\dots,0)$, so that $\left<\theta,X_r\right> $ has
the one dimensional density
$$
\int_{\R^{d-1}} p_r(x,y)\,dy \leq 
\frac{1}{b_r} \int_{|y| < r} p(x,y)\,dy \, \leq 2 \, \omega_{d-1} r^{d-1} M, \quad
x \in \R.
$$
Thus, for any unit vector $\theta$ in $\R^d$,
\be
M(\left<\theta,X_r\right>) \leq  2 \omega_{d-1} r^{d-1} M.
\en

Therefore, by Proposition 3.1 applied to the characteristic 
function $g_r(s)$ with its properties (4.3)-(4.4), it follows that
$$
1 - |g_r(s)|^2 \geq 
\frac{c}{\omega_{d-1}^2 r^{2(d-1)} M^2 \sigma^2}\,\min\{\sigma^2 s^2,1\}
$$
up to some absolute constant $c>0$. Using this in (4.2), we arrive at the
similar relation
$$
1 - |g(s)|^2 \geq 
\frac{c}{\omega_{d-1}^2 r^{2(d-1)} M^2 \sigma^2}\,\min\{\sigma^2 s^2,1\},
$$
that is,
$$
1 - |g(s)|^2 \geq 
\frac{c}{\omega_{d-1}^2 (2d)^{d-1} M^2\sigma^{2d}}\,
\min\{\sigma^2 s^2,1\}.
$$

To simplify the constants, recall the formula (2.3) and Batir's bounds 
for the Gamma function (\cite{Bat})
$$
\sqrt{2e}\,\Big(\frac{x}{e}\Big)^x \leq
\Gamma\Big(x + \frac{1}{2}\Big) \leq 
\sqrt{2\pi}\,\Big(\frac{x}{e}\Big)^x, \quad x \geq \frac{1}{2}.
$$
The lower bound gives
$$
\omega_{d-1} (2d)^{\frac{d-1}{2}} = 
\frac{\pi^{\frac{d-1}{2}}}{\Gamma(\frac{d+1}{2})}\, (2d)^{\frac{d-1}{2}}
 \leq \frac{1}{2\sqrt{\pi e d}}\,(4\pi e)^{d/2}.
$$
It remains to note that $1 - |g(s)|^2 \leq 2\,(1 - |g(s)|)$, and (4.1) follows.
\qed

\vskip10mm
\section{{\bf Maximum of Convolved Densities}}
\setcounter{equation}{0}

\vskip2mm
\noindent
Convolved densities are known to have improved smoothing properties. 
First, let us emphasize the following general fact.

\vskip5mm
{\bf Proposition 5.1.} {\sl If independent random vectors $X_1,\dots,X_m$
$(m \geq 2)$ with values in $\R^d$
have bounded densities, then the sum $S_m = X_1+ \dots + X_m$
has a bounded uniformly continuous density vanishing at infinity.
}

\vskip5mm
{\bf Proof.} Denote by $q_k$ the densities of $X_k$ and assume
that $q_k(x) \leq M_k$ for all $x \in \R^d$ with some constants
$M_k$ ($k \leq m$). By the Plancherel theorem, for the characteristic
functions $v_k(t) = \E\,e^{i\left<t,X_k\right>}$, we have
\bee
\int_{\R^d} |v_k(t)|^m\,dt
 & \leq &
\int_{\R^d} |v_k(t)|^2\,dt
 \, = \,
(2\pi)^d \int_{\R^d} q_k(x)^2\,dx \\
 & \leq &
(2\pi)^d \int_{\R^d} M_k q_k(x)\,dx \, = \, (2\pi)^d M_k,
\ene
where we used the property $|v_k(t)| \leq 1$, $t \in \R^d$.
Hence, by H\"older's inequality, the characteristic function
$f(t) = v_1(t) \dots v_m(t)$ of $S_m$ is integrable and has
$L^1$-norm
\begin{eqnarray}
\int_{\R^d} |f(t)|\,dt
 & \leq & 
\Big(\int_{\R^d} |v_1(t)|^m\,dt\Big)^{1/m} \dots
\Big(\int_{\R^d} |v_m(t)|^m\,dt\Big)^{1/m} \nonumber \\
 & \leq &
(2\pi)^d \, (M_1 \dots M_m)^{1/m} \, < \, \infty.
\end{eqnarray}
One may conclude that the random variable $S_m$ 
has a bounded, uniformly continuous density expressed by 
the inversion Fourier formula
\be
q(x) = \frac{1}{(2\pi)^d} \int_{\R^d} e^{-i\left<t,x\right>} f(t)\,dt,
\quad x \in \R.
\en
Since $f$ is integrable, it also follows that $q(x) \rightarrow 0$
as $|x| \rightarrow \infty$, by the Riemann-Lebesgue lemma.
\qed

\vskip2mm
Since, by (5.2),
$$
q(x) \leq \frac{1}{(2\pi)^d} \int_{\R^d} |f(t)|\,dt
$$
for all $x \in \R$, the inequality (5.1) also implies that
\be
M(S_m) \leq (M(X_1) \dots M(X_m))^{1/m}.
\en
However, the relation (5.3) 
does not correctly reflect  the behavior of $M(S_m)$ with respect 
to the growing parameter $m$, especially in the i.i.d. situation. 
A more precise statement from \cite{B-C} is described in the following 
relation, where the geometric mean of maxima is replaced with
the harmonic mean.

\vskip3mm
{\bf Proposition 5.2.} {\sl Given independent random vectors
$X_k$, $1 \leq k \leq m$, with values in $\R^d$, one has
\be
\frac{1}{M(S_m)^{2/d}} \geq \frac{1}{e} \sum_{k=1}^m 
\frac{1}{M(X_k)^{2/d}}.
\en
}

\vskip2mm
This bound may be viewed as a counterpart of the entropy power 
inequality in Information Theory. It is derived by applying
the Hausdorff-Young inequality with best constants (due to
Beckner and Lieb). The constant $1/e$ is optimal and is
attained asymptotically as $d \rightarrow \infty$ for random vectors
unifomly distributed on Euclidean balls. However, in dimension $d=1$,
it can be improved to $1/2$, which follows from results
due to Rogozin \cite{Ro} and Ball \cite{Ball1}.

One useful consequence of (5.4) is
the following bound on the $L^{2m}$-norms of characteristic functions.

\vskip2mm
{\bf Proposition 5.3.} {\sl If $f(t)$ is the characteristic function 
of a random vector $X$ in $\R^d$, then for any integer $m \geq 1$,
\be
\frac{1}{(2\pi)^d} \int_{\R^d} |f(t)|^{2m}\,dt \leq 
\Big(\frac{e}{2m}\Big)^{d/2}\, M(X).
\en
}

{\bf Proof.} We apply Proposition 5.2 to $2m$ summands 
$X_1,-X_1',\dots,X_m,-X_m'$, assuming that $X_k$, $X_k'$ are
independent copies of $X$. Introduce the symmetrized random vector
$\widetilde S_m = S_m - S_m'$, where $S_m'$ is an independent copy 
of $S_m$. By (4.4), we then get
$$
M(\widetilde S_m) \leq \Big(\frac{e}{2m}\Big)^{d/2}\, M(X).
$$
In addition, $\widetilde S_m$ has characteristic function $|f(t)|^{2m}$.
If $M(X)$ is finite, one may apply Proposition 5.1 and conclude 
that $\tilde S_m$ has a bounded continuous density $q_m(x)$ which
is vanishing at infinity. Moreover, $q_m(x)$ is maximized at $x = 0$, and 
its value at this point is described by the inversion formula (5.2) which gives
$$
M(\widetilde S_m) = q_m(0) = 
\frac{1}{(2\pi)^d} \int_{\R^d} |f(t)|^{2m}\,dt.
$$
\qed

\vskip2mm
Using (5.3), one can obtain a similar relation, but without
the factor $(\frac{e}{2m})^{d/2}$ in (5.5).

When $M(X)$ is finite and $m$ is large, the bound (5.5) may be 
considerably sharpened asymptotically with respect to $m$
when restricting the integration to the regions $|t| \geq \ep > 0$.

\vskip5mm
{\bf Proposition 5.4.} {\sl Let $f$ be the characteristic function 
of a random variable $X$ with covariance matrix $\sigma^2 I_d$ 
($\sigma>0$) and finite $M = M(X)$. For any $\ep > 0$ and $n \geq 2$,
\be
\int_{|t| \geq \ep} |f(t)|^n\,dt \leq 
\Big(\frac{8\pi^2 e}{n}\Big)^{d/2}\, M
\exp\Big\{-\frac{c^d n}{M^2\sigma^{2d}}\,\min(\sigma^2 \ep^2,1)\Big\}
\en
with some absolute constant $c>0$.
}

\vskip5mm
{\bf Proof.} Since the random vector $X$ has a density, we have
$\delta_f(\ep) = \max_{|t| \geq \ep} |f(t)| < 1$, 
by the Riemann-Lebesgue lemma. Moreover, by Proposition 4.1,
\be
\delta_f(\ep) \leq 1 - \frac{c^d}{M^2 \sigma^{2d}}\,\min(\sigma^2 \ep^2,1)
\en
with some absolute constant $c>0$. 
If $n=2$, we apply (5.5) with $m=1$ and use the property
that $M^2\sigma^{2d}$ is bounded away from zero, cf. (2.4). In the case 
$n \geq 3$, write $n = 2m + k$ with $k = 1$ or $k=2$ for $n \leq 5$ and 
with $m = [\frac{n}{3}]$, $k = n-2m$ for $n \geq 6$. Then, by (5.5) and (5.7),
\bee
\int_{|t| \geq \ep} |f(t)|^n\,dt 
 & \leq &
\delta_f(\ep)^k \int_{\R^d} |f(t)|^{2m}\,dt \\
  & \leq &
\Big(\frac{2\pi^2 e}{m}\Big)^{d/2}\, M
\exp\Big\{-\frac{c^d k}{M^2\sigma^{2d}}\,\min(\sigma^2 \ep^2,1)\Big\}.
\ene
It remains to note that $m \geq \frac{1}{4}\,n$ and $k\geq c_1 n$
for some absolute constant $c_1 > 0$.
\qed

\vskip10mm
\section{{\bf Normal Approximation for Products of Characteristic Functions}}
\setcounter{equation}{0}

\vskip2mm
\noindent
Let us now recall standard results about the approximation
of products of charactersitic functions. Consider the sum 
$$
S_n = \xi_1 + \dots + \xi_n
$$
of independent random variables $\xi_k$ with mean zero and 
standard deviations $b_k$ such that $b_1^2 + \dots + b_n^2 = 1$, 
so that $S_n$ has mean zero and variance one.

In terms of the characteristic functions $v_k(t) = \E\,e^{it \xi_k}$,
the characteristic function of $S_n$ represents the product
$$
f_n(t) = \E\,e^{it S_n} = v_1(t) \dots v_n(t), \quad t \in \R.
$$
Under higher order moment assumptions, it may be approximated 
by the standard normal characteristic function $g(t) = e^{-t^2/2}$ 
or its Edgeworth corrections on relatively large  $t$-intervals by 
means of the Lyapunov coefficients
$$
L_p = \sum_{k=1}^n \E\,|\xi_k|^p, \quad p > 2,
$$
provided that they are small. We only mention such results for 
the particular indexes $p = 3$ and $p = 4$.

\vskip5mm
{\bf Proposition 6.1.} {\sl With some absolute constants $C>0$ and $c>0$,
\be
|f_n(t) - e^{-t^2/2}|\leq C L_3\,|t|^3\,e^{-ct^2}, \quad |t| \leq \frac{1}{L_3}.
\en
Moreover, if $\E \xi_k^3 = 0$ for all $k \leq n$, then
\be
|f_n(t) - e^{-t^2/2}|\leq C L_4\,t^4\,e^{-ct^2}, \quad |t| \leq \frac{1}{\sqrt{L_4}}.
\en
}

\vskip2mm
The inequalities (6.1)-(6.2) are often stated in a slightly different form. 
For example, in Petrov (\cite{P}, p.\,109), the relation (6.1) is derived on 
a smaller interval $|t| \leq \frac{1}{4L_3}$ with $C = 16$ and $c = 1/3$.
As can be seen from the proof or properly modifying it, the interval of 
approximation can be increased to $|t| \leq \frac{c_0}{L_3}$ with some 
absolute constant $c_0>1$ at the expense of a smaller value of $c$ and 
a larger value of $C$. Similar relations with arbitrary real $p>2$ can be 
found in the review \cite{B1}.

In the non-interesting case, where $L_3$ or $L_4$ are greater than 1, 
(6.1)-(6.2) hold true on the larger interval $|t| \leq 1$. This can be seen 
from the Taylor integral fromula for the function $h(t) = f_n(t) - e^{-t^2/2}$ 
about the point $t=0$, using the property that the first 2 and 3 derivatives 
of $h(t)$ at the origin are respectively vanishing.

In general, the function $p \rightarrow L_p^{\frac{1}{p-2}}$ is 
non-decreasing with $L_p \geq n^{-\frac{p-2}{2}}$. In particular, 
\be
\frac{1}{\sqrt{n}} \leq L_3 \leq \sqrt{L_4}.
\en 
Since
$\E\,|\xi_k|^p \geq (\E\,\xi_k^2)^{p/2}$, we also have
\be
L_p \geq b_1^p + \dots + b_n^p \geq (\max_k b_k)^p.
\en

Now, if instead of $S_n$, one considers the normlized sum
$$
Z_n = \frac{X_1 + \dots + X_n}{\sqrt{n}},
$$
these results can be correspondingly reformulated. Then one should assume
that the independent random variables $X_k$ have mean zero and standard 
deviations $\sigma_k$ such that $\sigma_1^2 + \dots + \sigma_n^2 = n$.
In this setting, it is more natural to represent the Lyapunov coefficients as
$$
L_p = n^{-\frac{p-2}{2}}\,\beta_p, \quad 
\beta_p = \frac{1}{n} \sum_{k=1}^n \E\,|X_k|^p.
$$
Thus, $\beta_p = \E\,|X|^p$ when $X_k$ are independent copies of
a random variable $X$.

In general, from (6.3)-(6.4) it follows that
\be
1 \leq \beta_3 \leq \sqrt{\beta_4}, \quad 
1 \leq \max_k\, \sigma_k \leq (n \beta_p)^{1/p}.
\en

Finally, let us state Proposition 6.1 once more about the normalized sums.

\vskip5mm
{\bf Proposition 6.2.} {\sl With some absolute constants $C>0$ and $c>0$,
the characteristic function $f_n(t)$ of $Z_n$ satisfies
\be
|f_n(t) - e^{-t^2/2}|\leq \frac{C \beta_3}{\sqrt{n}}\,|t|^3\,e^{-ct^2}, \quad 
|t| \leq \frac{\sqrt{n}}{\beta_3}.
\en
Moreover, if $\E X_k^3 = 0$ for all $k \leq n$, then
\be
|f_n(t) - e^{-t^2/2}|\leq \frac{C \beta_4}{n}\,t^4\,e^{-ct^2}, \quad
|t| \leq \frac{\sqrt{n}}{\sqrt{\beta_4}}.
\en
}

\vskip5mm
\section{{\bf Refinement of Theorem 1.1}}
\setcounter{equation}{0}

\vskip2mm
\noindent
Let us return to the setting of Theorem 1.1 and assume that the independent
random vectors $X_k$ have mean zero, covariance matrix $\sigma_k^2 I_d$,
and finite maxima of densities $M_k = M(X_k)$. Recall the notations
$$
\beta_3 = \sup_{|\theta| = 1}\,
\bigg[\,\frac{1}{n} \sum_{k=1}^{n} \E\,|\left<\theta,X_k\right>|^3\bigg], 
\quad
\beta_4 = \sup_{|\theta| = 1}\,
\bigg[\,\frac{1}{n} \sum_{k=1}^{n} \E \left<\theta,X_k\right>^4\bigg],
$$
and
$$
\Delta_n = \sup_x |p_n(x) - \varphi(x)|.
$$
Here we prove this theorem in a somewhat more general form (although more complicated).

\vskip5mm
{\bf Theorem 7.1.} {\sl With some positive absolute constants $C,c$, the density
$p_n$ of $Z_n$ satisfies
\be
\Delta_n \leq \frac{C^d\beta_3}{\sqrt{n}} + C^d\, (M_1 \dots M_n)^{\frac{1}{n}}
\exp\Big\{-c^d\,\sum_{k=1}^n 
\frac{1}{M_k^2 \sigma_k^{2d}}\,\min(\sigma_k^2/\beta_3^2,1)\Big\}.
\en
Moreover, if $\E \left<\theta,X_k\right>^3 = 0$ for all $\theta \in \R^d$ and 
$k \leq n$ (in particular, if the distributions of $X_k$ are symmetric), then
\be
\Delta_n \leq \frac{C^d\beta_4}{n} + C^d\, (M_1 \dots M_n)^{\frac{1}{n}}
\exp\Big\{-c^d\,\sum_{k=1}^n 
\frac{1}{M_k^2 \sigma_k^{2d}}\,\min(\sigma_k^2/\beta_4,1)\Big\}.
\en
}

\vskip2mm
{\bf Proof.} Since necessarily $M \geq c^d$ and 
$\sqrt{\beta_4} \geq \beta_3 \geq 1$ (cf. (2.4)),
the inequalities (7.1)-(7.2) are fulfilled automatically for $n=1$. 
So, assume that $n \geq 2$.

In terms of the characteristic functions $v_k(t) = \E\,e^{i\left<t,X_k\right>}$,
the characteristic function of $Z_n$ is given by the product
$$
f_n(t) = \E\,e^{i\left<t,Z_n\right>} =
v_1\Big(\frac{t}{\sqrt{n}}\Big) \dots v_n\Big(\frac{t}{\sqrt{n}}\Big).
$$
Applying H\"older's inequality and then Proposition 5.3, we get
$$
\int_{\R^d} |f_n(t)|\,dt \leq \prod_{k=1}^n
\bigg(\int_{\R^d} 
\Big|v_k\Big(\frac{t}{\sqrt{n}}\Big)\Big|^n\,dt\bigg)^{1/n} < \infty.
$$
Hence, one may apply the Fourier inversion
formula to represent the densities of $Z_n$ as
$$
p_n(x) = \frac{1}{2\pi} \int_{\R^d} e^{-i\left<t,x\right>} f_n(t)\,dt, \quad
x \in \R^d.
$$
Using a similar representation for the standard normal density, we get
\be
\Delta_n \leq 
\frac{1}{2\pi} \int_{\R^d} |f_n(t) - e^{-|t|^2/2}|\,dt.
\en

Next, we split integration in (7.3) into the two regions.
From Proposition 5.4 applied to the characteristic function $v_k(t)$ of 
$X_k$, we obtain that, for any $\ep>0$,
$$
\int_{|t| \geq \ep} |v_k(t)|^n\,dt \leq 
\Big(\frac{C_0}{\sqrt{n}}\Big)^d\, M_k
\exp\Big\{-c^d\, \frac{n}{M_k^2 \sigma_k^{2d}}\,\min(\sigma_k^2 \ep^2,1)\Big\},
$$
where $C_0 = 2\pi\sqrt{2e}$. Equivalently,
$$
\int_{|t| \geq \ep\sqrt{n}} \Big|v_k\Big(\frac{t}{\sqrt{n}}\Big)\Big|^n\,dt \leq 
C_0^d\, M_k
\exp\Big\{-c^d\, \frac{n}{M_k^2 \sigma_k^{2d}}\,\min(\sigma_k^2 \ep^2,1)\Big\}.
$$
Applying once more H\"older's inequality
$$
\int_{|t| \geq \ep\sqrt{n}} |f_n(t)|\,dt \leq \prod_{k=1}^n
\bigg(\int_{|t| \geq \ep\sqrt{n}} 
\Big|v_k\Big(\frac{t}{\sqrt{n}}\Big)\Big|^n\,dt\bigg)^{1/n},
$$
we then get
$$
\int_{|t| \geq \ep\sqrt{n}} |f_n(t)|\,dt \, \leq \,
C_0^d\, (M_1 \dots M_n)^{\frac{1}{n}} \exp\Big\{-c^d\,\sum_{k=1}^n 
\frac{\min(\sigma_k^2 \ep^2,1)}{M_k^2 \sigma_k^{2d}}\Big\}.
$$
Since
$$
\int_{|t| \geq \ep\sqrt{n}} e^{-|t|^2/2}\,dt \, \leq \, C^d e^{-n\ep^2/4}
$$
for some absolute constant $C>0$, it follows that
\begin{eqnarray}
\int_{|t| \geq \ep\sqrt{n}} \big|f_n(t) - e^{-|t|^2/2}\big|\,dt 
 & \leq &
C^d e^{-n\ep^2/4} \nonumber \\
 & & \hskip-20mm + \ 
C^d\, (M_1 \dots M_n)^{\frac{1}{n}} \exp\Big\{-c^d\,\sum_{k=1}^n 
\frac{\min(\sigma_k^2 \ep^2,1)}{M_k^2 \sigma_k^{2d}}\Big\}.
\end{eqnarray}

We now turn to the other region $|t| < \ep \sqrt{n}$.
By the mean zero and isotropy assumption on the distribution of $X_k$, 
we have $\E \left<\theta,X_k\right> = 0$ and
$\E \left<\theta,X_k\right>^2 = 1$ for any unit vector $\theta$
in $\R^d$. Therefore, we are in position to apply Proposition 6.2 to 
the random variables $\left<\theta,X_k\right>$. Write $t = s\theta$ 
for $s>0$ and $|\theta| = 1$.
Since $s \rightarrow f_n(s\theta)$ represents the characteristic function
of the normalized sum of $\left<\theta,X_k\right>$, it satisfies, by (6.6),
\be
|f_n(s\theta) - e^{-s^2/2}| \leq 
\frac{C \beta_3(\theta)}{\sqrt{n}}\,s^3\,e^{-cs^2}, \quad 
s \leq \frac{\sqrt{n}}{\beta_3(\theta)},
\en
where 
$$
\beta_3(\theta) = \frac{1}{n} \sum_{k=1}^n \E\,|\left<\theta,X_k\right>|^3.
$$
In (7.5), $\beta_3(\theta)$ can be replaced with the larger value $\beta_3$, 
which leads to
$$
|f_n(t) - e^{-|t|^2/2}| \leq \frac{C \beta_3}{\sqrt{n}}\,|t|^3\,e^{-c|t|^2}, \quad 
|t| \leq T_n = \frac{\sqrt{n}}{\beta_3}.
$$
This readily yields
$$
\int_{|t| \leq T_n} |f_n(t) - e^{-|t|^2/2}|\,dt \leq 
\frac{C^d\beta_3}{\sqrt{n}}.
$$
with some absolute constant $C$. One can now combine this inequality 
with (7.4) by choosing $\ep = \frac{1}{\beta_3}$. Since $\beta_3 \geq 1$, 
the resulting inequality in (7.3) yields (7.1).

In the second scenario, we similarly apply the inequality (6.7)
and combine it once more with (7.4) by choosing $\ep = \frac{1}{\sqrt{\beta_4}}$. 
\qed

\vskip7mm
\section{{\bf Proof of Theorem 1.1 and Corollary 1.2}}
\setcounter{equation}{0}

\vskip2mm
\noindent
The right-hand sides in (7.1)-(7.2) may be bounded and simplified in terms 
of the functionals
$$
M = \max_k M(X_k), \quad \sigma^2 = \max_k \sigma_k^2.
$$

\vskip2mm
{\bf Proof of Theorem 1.1.}
The $k$-term in the sum in (7.1) represents a decreasing function with
respect to $\sigma_k^2$ and $M_k^2$. Hence the second summand on the
right-hand side in (7.1) does not exceed
$$
C^d M
\exp\Big\{-\frac{c^d n}{M^2 \sigma^{2d}}\,\min(\sigma^2/\beta_3^2,1)\Big\}.
$$
Moreover, since $\sigma^2 \geq 1$ and $\beta_3 \geq 1$, necessarily
$\min(\sigma^2/\beta_3^2,1) \geq \beta_3^{-2}$. This further simplifies
the latter expression to
$$
C^d M
\exp\Big\{-\frac{c^d n}{M^2 \sigma^{2d} \beta_3^2}\Big\}.
$$
Using $e^{-x} < x^{-1/2}$ ($x>0$), from (7.1) we obtain that
$$
\Delta_n \leq \frac{C^d\beta_3}{\sqrt{n}} + C^d\,
\frac{M^2 \sigma^d \beta_3}{\sqrt{n}}.
$$
Here, the first term on the right-hand side is dominated by 
the second one up to the multiple, so that
the above estimate is simplified to (1.1).

With similar arguments, the second summand on the right-hand side
in (7.2) does not exceed
$$
C^d M
\exp\Big\{-\frac{c^d n}{M^2 \sigma^{2d} \beta_4}\Big\}.
$$
Using $e^{-x} < x^{-1}$ ($x>0$), from (7.2) we therefore obtain that
$$
\Delta_n \leq \frac{C^d\beta_4}{n} + C^d\, \frac{M^3 \sigma^{2d} \beta_4}{n}.
$$
Again, the first term on the right-hand side is dominated by 
the second one up to the multiple, so that
the above estimate is simplified to (1.2). This proves Theorem 1.1.
\qed

\vskip5mm
{\bf Proof of Corollary 1.2.}
First we need to mention that, if a random vector $X$
has a log-concave density, so does any linear functional $\xi = \left<\theta,X\right>$.
This is a consequence of a well-known characterization of log-concave measures 
due to Borell \cite{Bor}. He also derived a large deviation bound for norms 
under log-concave measures, which implies in dimension one that
$L^p$-norms of $\xi$ are equivalent to each other. More precisely,
$$
(\E\,|\xi|^p)^{1/p} \leq Cp\,(\E \xi^2)^{1/2}
$$ 
for all $p > 2$ with some absolute constant $C$.
Applying this with $X = X_k$ in the case $\sigma_k = 1$, we conclude that 
$\beta_3$ and $\beta_4$ are bounded by absolute constants. Hence,
the inequalities (1.1)-(1.2) are respectively simplified to
$$
\Delta_n \leq \frac{C^d M^2}{\sqrt{n}}, \quad
\Delta_n \leq \frac{C^d M^3}{n}.
$$
To prove (1.3)-(1.4), it remains to recall the bound (2.5) which gives
$M \leq K_d^d$, where $K_d$ is the maximal isotropic constant for 
the class of log-concave probability distributions on $\R^d$.

\end{document}